\documentclass{amsart}

\usepackage{bbm}
\usepackage{amssymb}
\usepackage{tikz,ifthen}
\usepackage{url}
\usepackage{units}
\usepackage{xspace}
\usepackage{paralist}
\usepackage{cite}
\usepackage[colorlinks]{hyperref}
\usepackage{breakurl}
\usepackage{listings}
\newcommand{\ve}{\boldsymbol{e}}

\newcommand{\vv}{\boldsymbol{v}}

\newcommand{\R}{\mathbbm R}

\newcommand{\Z}{\mathbbm Z}
\newcommand{\N}{\mathbbm N}

\DeclareMathOperator{\aff}{aff}

\DeclareMathOperator{\lvol}{lvol}
\DeclareMathOperator{\pyr}{pyr}
\DeclareMathOperator{\conv}{conv}
\DeclareMathOperator{\faces}{F}
\newcommand{\polymake}{\texttt{polymake}\xspace}
\DeclareMathOperator{\ehr}{ehr}

\begin{document}

\title[Defect Polytopes and Counter-Examples With polymake]{Defect Polytopes and Counter-Examples\\ With polymake}

\author{Michael Joswig}
\thanks{Michael Joswig is supported by DFG Research Group ``Polyhedral Surfaces'' (FG 565) and DFG Priority Program ``Algorithmic and Experimental Methods'' (SPP 1489)}

\author{Andreas Paffenholz}
\thanks{Andreas Paffenholz is supported by the DFG Priority Program ``Algorithmic and Experimental Methods'' (SPP 1489)}

\address{TU Darmstadt, Dolivostr. 15, 64293 Darmstadt, Germany}
\email{\{joswig,paffenholz\}@mathematik.tu-darmstadt.de}

\maketitle

\begin{abstract}
  It is demonstrated how the software system \polymake can be used for computations in toric
  geometry.  More precisely, counter-examples to conjectures related to A-determinants and defect
  polytopes are constructed.
\end{abstract}

\section{Defect Polytopes}

The purpose of this note is to demonstrate how \polymake's features can be used to establish
explicit counter-examples to some conjectures on defect polytopes in toric geometry.  More
precisely, we consider several invariants of lattice polytopes proposed by Di
Rocco~\cite{1098.14039}, show a way to compute them with \polymake, and use this to give answers to
questions posed in Section~6 of that paper.  For basic properties of polytopes and related toric
varieties we refer to Ewald~\cite{Ewald}.

We start by introducing some notation. Let $P$ be a lattice polytope, i.e., a polytope whose
vertices are contained in $\Z^d$. Let $\faces(P)$ the set of all
non-empty faces of $P$, including $P$ itself.  Further, we let $\faces_P(k)$ be the set of
$k$-dimensional faces of $P$. For $t\in \N$ we define
\begin{equation}\label{eq:c}
  c_t(P) \ := \ \sum_{k=0}^d (-1)^{d-k}\frac{(k+t)!}{k!}\sum_{F\in\faces_P(k)}\lvol(F) \, ,
\end{equation}
where $\lvol(F)$ is the normalized volume of $F$ in the lattice
\[
  \Z^F \ := \ \Z^d\cap\aff(F) \, ;
\]
see \cite[p.~101]{1098.14039}. For smooth polytopes $P$ the invariant $c_1(P)$ is a function that
records the degree of homogeneity of a certain rational function, the \emph{$A$-determinant}, where
$A:=P\cap \Z^d$ is the set of lattice points in $P$. Starting with the seminal
monograph~\cite{MR2394437} this invariant was extensively studied; see
also~\cite{0810.4996,0807.3163,1001.2792,Curran2007115} for recent results, and generalizations to
the singular case.  By \cite[\S11, Thm.~1.6]{MR2394437} the number $c_1(P)$ is non-negative for
simple lattice polytopes whose associated toric variety is quasi-smooth. (and thus, in particular,
for smooth polytopes). By \cite[Pro.~3.2]{1001.2792} this also holds for arbitrary lattice
simplices, while it is open for more general simple polytopes.  A polytope $P$ is a \emph{defect
  polytope}, if $P$ is a smooth lattice polytope and $c_1(P)=0$.  In terms of the associated toric
variety $P$ is a defect polytope if the dual variety is not a hypersurface.  By a result of
Dickenstein and Nill~\cite[Thm.~1.6]{1098.14039} any defect polytope is a smooth strict Cayley
polytope. A \emph{strict Cayley polytope} is a polytope that is affinely isomorphic to
\begin{align*}
  Q_0\star\cdots\star Q_k=\conv(Q_0\times\{\ve_0\},\ldots, Q_k\times\{\ve_k\})
\end{align*}
where $\ve_j$ is a lattice basis of $\Z^{k+1}$ and $Q_0, \ldots, Q_k$ are strictly isomorphic
lattice polytopes in $\R^m$ (i.e.\ having the same normal fan) such that $\dim(\aff(Q_0,\ldots,
Q_k))=m$.

\section{The Software}

\polymake is a software system for computations in geometric combinatorics and related areas. The
project was initiated in 1995 by Gawrilow and the first author~\cite{DMV:polymake}, and many people
helped to continuously expand it since. Recently, two important additions to the system have been
accomplished~\cite{FPSAC2009}:
\begin{compactenum}
\item \polymake now comes with an interactive shell similar to most computer algebra systems.
\item \polymake has been extended to allow computations specific to the class of \emph{lattice
    polytopes}, i.e., convex polytopes with integral vertex coordinates, and their relation to
  combinatorial commutative algebra, toric geometry, and integer programming.
\end{compactenum}
The latest release 2.9.9 of \polymake was published on November 9, 2010 and can be obtained from
\url{http://www.polymake.org}. It is distributed as source code and precompiled binaries for Mac
OS~X and several Linux distributions. It is released under the GNU GPL, version 3.

\polymake's functionality is organized in various \textit{applications}.  Currently, these are
\texttt{polytope} for computations with convex polyhedra, \texttt{matroid} and \texttt{graph} for
purposes revealed by their names, \texttt{topaz} (which is short for ``\underline{top}ology
\underline{a}pplication \underline{z}oo'') for finite simplicial complexes, and \texttt{tropical}
for tropical geometry. Each application centers around \textit{objects}, a representation for the
basic mathematical objects dealt with, e.g., \texttt{Polytope<Rational>} or \texttt{LatticePolytope}
in the application \texttt{polytope}.  Technically, objects are organized in a class hierarchy
written in (a slightly extended dialect of) Perl.  These extensions include the possibility to use
C++-style template parameters (e.g., \texttt{<Rational>} to specify a coordinate domain for
polytopes) and a shared-memory communication model to interface with compiled C++-code.
Semantically, objects are defined by their properties with can be used to derive further properties
according to an extendible set of rules.  For more about the general ideas behind \polymake the
reader should consult~\cite{DMV:polymake}, while \cite{FPSAC2009} is a more recent account with a
focus on lattice polytopes.

\section{Computations and Code Fragments}

We prepare a function, written in \polymake's Perl dialect, that computes the function $c$
from~\eqref{eq:c}. It takes a polytope $P$ in the first, and the parameter $t$ in the second
argument.
\begin{footnotesize}
\begin{verbatim}
sub ct_invariant {
 my ($P, $t) = @_;
 my $v = $P->VERTICES;
 my $hd = $P->HASSE_DIAGRAM;
 my $sign = 1;
 my $c = new Integer(0);

 for (my $d = $P->DIM; $d > 0; --$d) {
  foreach (@{$hd->nodes_of_dim($d)}) {
   my $F = new Polytope(VERTICES=>
              $v->minor($hd->FACES->[$_],All));
   my $vol = $F->LATTICE_VOLUME;
   $c += $sign*fac($d+$t)/fac($d)*$vol;
  }

  $sign = -$sign;
 }

 $c += $sign*fac($t)*$P->N_VERTICES;
 return $c;
}
\end{verbatim}
\end{footnotesize}
The code of all functions used in this section, together with some explanations how to use them, is
also available as a tutorial in the \polymake wiki at
\url{http://www.polymake.org/doku.php/tutorial/a_determinants}.  \polymake\ has interfaces to
several other software packages for the computation of some special properties, such as the volume
with respect to a given lattice, in the example above.  \polymake 2.9.9 uses
\texttt{Normaliz}~\cite{normaliz2} for this type of computation by default.

Now we can interactively apply the function \texttt{ct\_invariant} to some examples, e.g.,
\cite[Example 2, p.~86]{1098.14039}. Each input line to the \polymake shell is preceded with the
name of the currently active application, which is \texttt{polytope} throughout this extended
abstract:
\begin{footnotesize}
\begin{verbatim}
polytope > $S = simplex(2);
polytope > $P = prism($S);
polytope > print $P->SMOOTH;
1
polytope > print ct_invariant($P,1);
0
\end{verbatim}
\end{footnotesize}
Here we have assigned $\conv\{(0,0),(1,0),(0,1)\}$, the standard simplex in $\R^2$, to the variable
\texttt{\$S} and defined \texttt{\$P} to be the prism over \texttt{\$S} (with distance $1$ between
base and top). We can check for non-singularity of the associated toric variety by asking for the
property \texttt{SMOOTH} in \polymake.  The output ``1'' represents the boolean value
``true''. Calling the function \texttt{ct\_invariant} defined above shows that, indeed,
$c_1(P)=0$.  So the prism over the standard triangle is a defect polytope. It is a Cayley polytope
of $3$ segments. We can repeat the same computation for the hypersimplex $\Delta(3,6)$, which is not
simple and hence not smooth:
\begin{footnotesize}
\begin{verbatim}
polytope > print ct_invariant(hypersimplex(3,6),1);
136
\end{verbatim}
\end{footnotesize}
Notice that the value given in \cite[Ex.~10]{1098.14039} is not correct.

Using the connection between lattice polytopes and toric varieties Di Rocco proves that $c_t(P)\ge
0$ for $t\ge 1$ and any smooth lattice polytope $P$~\cite[Cor.\ 4]{1098.14039}.

\section{Counter-Examples}

In the sequel we use \polymake to give answers to three questions
raised in \cite{1098.14039}.

\subsection{Lattice polytopes with $c_1$ negative}
Conjecture~4 in that paper asks if $c_1(P)\ge 0$ for any lattice
polytope, not necessarily smooth.  We introduce some more notation in
order to give counter-examples.  Let $P$ be a lattice polytope in
$\R^d$. The \emph{lattice pyramid} of $P$ is the lattice polytope
\begin{align*}
  \pyr(P)\ := \ \conv( P\times \{0\}\; \cup \; \vv\times\{1\} )\,,
\end{align*}
where $\vv$ is a vertex of $P$ (any other lattice point at distance $1$ from $\R^d\times \{0\}\subseteq \R^{d+1}$ gives an equivalent lattice polytope). We define the \emph{$r$-fold lattice   pyramid} of $P$ inductively via
\begin{align*}
  \pyr^r(P)\; &:=\; \pyr(\pyr^{r-1}(P))&&\text{for $r\ge 1$,}\\
\pyr^{0}(P)\; &:=\; P\;.
\end{align*}
We define another small function to compute $r$-fold pyramids.
\begin{footnotesize}
\begin{verbatim}
sub r_fold_pyr {
 my ($P, $r) = @_;
 for (; $r>0; --$r) {
  $P = lattice_pyramid($P);
 }
 return $P;
}
\end{verbatim}
\end{footnotesize}
We can apply this to the standard unit cube, and compute some of the invariants $c_t$ for lattice
pyramids over the cube.
\begin{footnotesize}
\begin{verbatim}
polytope > $C=cube(3,0);
polytope > print ct_invariant($C,0);
-2
polytope > print ct_invariant($C,1);
4
polytope > print ct_invariant(r_fold_pyr($C,1),1);
-1
polytope > print ct_invariant(r_fold_pyr($C,2),2);
-2
polytope > print ct_invariant(r_fold_pyr($C,3),3);
-6
\end{verbatim}
\end{footnotesize}
With the knowledge of this example it is not hard to see that 
\begin{align*} 
  c_r(\pyr^r(P))\ & =\ r!c_0(P)+(-1)^{d+1}r!\;.
\end{align*}
for any $d$-dimensional lattice polytope $P$ and $r>0$. We have seen in the previous computation
that $c_0(C)=-2$ for the standard unit cube $C$, so that
\begin{align*}
  c_r(\pyr^r(C))\ =\ r!c_0(C)+r! \ = \ -r! \ < \ 0\,.
\end{align*}

\subsection{Convolutions of Ehrhart polynomials}
Conjecture~5 in \cite{1098.14039} asks whether all coefficients of the
polynomial
\begin{align*}
  f(P,t) \ &:= \ \sum_{k=0}^d(-t)^{(d-k)}(k+1)!\sum_{F\in \faces_P(k)}|Z^F\cap tF|\\
  &= \ \sum_{k=0}^d\sum_{F\in \faces(P)}(k+1)!\ehr(F,t)t^{(d-k)}
\end{align*}
in the indeterminate~$t$ are non-negative. Here $\ehr(F,t)$ is the Ehrhart polynomial of $F$, and
the coefficients of $f(P,t)$ are integral, as $k!\ehr(F,k)$ has integral coefficients for any
$k$-dimensional lattice polytope.  Again we will provide counter-examples.  In \polymake, we can use
the following function to compute the coefficients of $f(P,t)$.
\begin{footnotesize}
\begin{verbatim}
sub f_poly_coeff {
 my ($P) = @_;
 my $d = $P->DIM;
 my $v = $P->VERTICES;
 my $hd = $P->HASSE_DIAGRAM;
 my $sign = 1;
 my $f = new Vector<Integer>($d+1);

 for (my $k = $d;  $k > 0; --$k) {
  foreach (@{$hd->nodes_of_dim($k)}) {
   my $q = new Polytope(VERTICES=>
              $v->minor($hd->FACES->[$_],All));
   my $h = (zero_vector<Rational>($d-$k))
                 |$q->EHRHART_POLYNOMIAL_COEFF;
   $f += $sign*fac($k+1)*$h;
  }
  $sign = -$sign;
 }

 $f += $sign*$P->N_VERTICES
                 *unit_vector<Rational>($d+1,$d);
 return $f; 
}
\end{verbatim}
\end{footnotesize}
The returned vector lists the coefficients of the polynomial with increasing degree, i.e., the last
entry is the leading coefficient of the polynomial.

We can use the same example as before to show that the coefficients need not be non-negative.  We
start out with the $3$-dimensional unit cube stored in \texttt{\$C}.
\begin{footnotesize}
\begin{verbatim}
polytope > print f_poly_coeff($C);
24 36 24 4
polytope > print f_poly(lattice_pyramid($C));
120 192 114 32 -1
\end{verbatim}
\end{footnotesize}
It is not hard to see that the leading coefficient of this polynomial coincides with $c_1(P)$.
Looking at iterated pyramids, by re-using our function \texttt{r\_fold\_pyr}, we obtain:
\begin{footnotesize}
\begin{verbatim}
polytope > print f_poly_coeff(r_fold_pyr($C,3));
5040 9060 5538 1698 188 -3 0
polytope > print f_poly_coeff(r_fold_pyr($C,5));
362880 717696 491304 163056 28086 1490 -15 0 0
\end{verbatim}
\end{footnotesize}
We are not aware of an example where another but the leading coefficient of the polynomial is
negative.

\subsection{Cayley polytopes}
It is finally asked in \cite[p.~103]{1098.14039} whether $c_1(P)=0$ implies that $P$ is a strict
Cayley polytope as defined above. This is not true, as the example of a two-fold pyramid over the
unit square shows:
\begin{footnotesize}
\begin{verbatim}
polytope > $Q=cube(2,0);
polytope > print ct_invariant(r_fold_pyr($Q,2),1);
0
\end{verbatim}
\end{footnotesize}

However, in the above definition for a Cayley polytope we could drop the condition that all factors
have the same normal fan. With this slightly more general notion the $2$-fold lattice pyramid over
the $0/1$-square is a non-strict Cayley polytope of two segments and two points.  On the other hand,
it is known~\cite[Prop.~2.14]{0810.4996} that the set of lattice points of a defect polytope is
contained in two parallel hyperplanes with distance one, so any defect polytope is always a
non-strict Cayley polytope over a segment.

%

\goodbreak


\providecommand{\bysame}{\leavevmode\hbox to3em{\hrulefill}\thinspace}
\providecommand{\MR}{\relax\ifhmode\unskip\space\fi MR }
\providecommand{\MRhref}[2]{%
  \href{http://www.ams.org/mathscinet-getitem?mr=#1}{#2}
}
\providecommand{\href}[2]{#2}

\end{document}